\documentclass{article}
\usepackage{graphicx} 
\usepackage{amsmath,amsthm,amssymb}
\usepackage{mathtools}
\usepackage{abstract}

\newtheorem{theorem}{Theorem}
\newtheorem{lemma}[theorem]{Lemma}

\def\cC{{\cal C}}

\begin{document}

\title{Complexity of graph evolutions \thanks{preliminary version, November 29, 2024}}
\author{Jeffrey P. Gao \& Paul C. Kainen\\  {\tt \{jpg135, kainen\}@$\,$georgetown.edu}}
\date{}

\maketitle

\noindent
\begin{abstract} \noindent
A permutation of the elements of a graph is a {\it construction sequence} if no edge is listed before either of its endpoints. The complexity of such a sequence is investigated by finding the delay in placing the edges, an {\it opportunity cost} for the construction sequence.   Maximum and minimum cost c-sequences are provided for a variety of graphs and are used to measure the complexity of graph-building programs.
\end{abstract}

\section{Introduction}

We study the ways in which a graph can {\it evolve} from a set of vertices and edges in such a way that, at every stage, one has a graph.  Such an evolution is called ``a construction sequence''  in \cite{kp} which counts these sequences.  Using a sequence of consecutive integers to enumerate the stages allows combinatorial measures.

The concept of a graph is enriched by this ontological model:  A graph comes into being as a sequence of elementary events when the elements of a graph, its vertices and edges, are appended to a list in which the initial subsequences always describe a graph.  Because the number of such construction sequences grows very rapidly with the graph, it may be useful to think probabilistically

An evolving graph is a process that sequentially chooses vertices and edges from a given graph, or from a distribution of possible graphs, subject only to the constraint that each new edge is attached to two existing vertices; in the construction sequence, {\it edges follow both their endpoints}.  Thus, our definition of an evolving graph is a linear order on the set of vertices and edges extending the  partial order of containment of vertices in edges of which they are endpoints.

Suppose a star-graph $K_{1,n}$ is evolving in the plane. One might draw different conclusions if the order of appearance of the $n$ edges (the {\it spokes}) were clockwise rather than random.
Thus, in choosing a c-sequence, the process is giving more information than just the resulting graph.  

Even without having external data regarding the graph's layout in space (i.e., when and where edges are placed), one can record an {\it internal} parameter: the delay between when each edge is placed and when it {\it could} $\,$have been placed (i.e., one step after  the later of  its endpoints has been placed).  One can then combine the internal parameters by summation to obtain the total global delay
of the construction sequence which we view as its {\it cost}.

The idea of counting the number of linear extensions of a partial order has been developed by Stanley \cite{stan} in a very general setting but he only evaluated the number of c-sequences for paths in \cite{stan2}, where it is shown that the {\it construction number} of a path with $n$ vertices is the $n$-th tangent number \cite[A000182]{oeis}.  For the construction numbers of paths, stars, complete graphs, et al., see \cite{kp, monthly}.


Here, we study the complexity of c-sequences by considering the total delay.  Edges, like the verbs in a German sentence, can be placed at the end of the list. In the extreme case of what we call an ``easy'' sequence, all the vertices are listed first, followed by the edges.  Such an evolution is required in order to achieve the greatest possible delay
in placing the edges. 

An edge is {\it available} for a c-sequence if both its endpoints have appeared.
To get the least total delay, all edges that are available should be placed as soon as possible; this leads to what we call ``greedy'' c-sequences which must be followed by the evolution to achieve the least delay.  Thus, in a conceptually efficient logical development, after two statements have been made, any implications between them should be proved, before other unrelated statements are given.

For a connected graph $G$, one may consider  ``nearly connected'' construction sequences \cite{kp}, where the evolving sequence of subgraphs has all members with at most two connected components. In a number of examples, we can show that minimum cost c-sequences for connected graphs are nearly connected as well as greedy.
However, there are minimum cost c-sequences for some graphs that are not nearly connected.  
As in the games of Go and Hex, one can play isolated stones advantageously before connecting them up.

In this paper, we restrict to discrete time. With a {\it continuous} model for time, R. Stong \cite{ks} gave an elegant proof that among $n$-vertex trees, the path and the star have smallest and largest construction number.  We show below that for the max cost function on $n$-vertex trees, path and star are similarly extremal.

Cost as a model for complexity could be applied to construction sequences for hypergraphs or cellular complexes, and we believe this idea may be helpful in finding efficient deployments within operational research.




    Section 2 has basic definitions, several fundamental lemmas, and examples. In Section 3,  we find maximum cost c-sequences for various regular and irregular graphs.  Minimum cost c-sequences are treated in Section 4, and we conclude in Section 5 with a discussion of potential applications and alterations to our assumptions. 
A brief appendix summarizes results for a number of other graphs.

\section{Definitions, lemmas, and examples}

For a finite set $A$, we write $|A|$ for the cardinality. Let $[n] := \{1,2, \ldots, n\}$.  
If $G = (V,E)$ is a graph, we let $p := |V|$ and $q := |E|$.  Put $\ell := p+q$ and call $V \sqcup E$ the set of {\bf elements} of $G$. For undefined graph terminology, see e.g. \cite{harary}.  We write $\odot$ for concatenation of strings.

A {\bf construction sequence} for $G = (V,E)$ is a bijection $x: [\ell] \to V \sqcup E$ such that if $e = uw$ is any edge of $G$ and $u$ and $w$ are its endpoints, then
\begin{equation}
x^{-1}(e) \,>\, \max \{x^{-1}(u), x^{-1}(w)\}. 
\label{eq:cs}
\end{equation}
If $t$ is an element of $G$ and $x$ is a construction sequence (c-sequence) for $G$ with $x^{-1}(t) = i$, we say that $t$ is {\bf placed} at time $i$.  Condition (\ref{eq:cs}) says that an edge must be placed after both of its endpoints.
The {\bf construction number} $\cC(G)$ of $G$ is $\cC(G) :=|\cC(G)|$, where
$\cC(G)$ is the set of all c-sequences for $G$. 

For $x \in \cC(G)$, let $x_{(i)} := (x_1, \ldots, x_i)$ be an initial subsequence of $x$. When $1 \leq i \leq \ell$, we write $G_i := G(\{(x_1, \ldots, x_i\}) \subseteq G$ for the {\bf initial subgraph} of $x$  induced by $\{x_1,\ldots, x_i\}$. The initial subsequence $x_{(i)}$ is a c-sequence for $G_i$.

Using (\ref{eq:cs}) define the {\bf cost of an edge $e = uw$ in a c-sequence $x$}, 
\begin{equation}
\nu(e, x) = (x^{-1}e - x^{-1}u) + (x^{-1}e - x^{-1}w).
\end{equation}
The {\bf cost of the c-sequence} $x$ for $G$ is then given by summing over the edges
\begin{equation}
\nu(x) = \sum_{e \in EG} \nu(e, x).
\end{equation}
Taking maximum and minimum over all construction sequences, one gets the {\bf max and min cost to build a graph} $G$:
\begin{equation}
\nu^*(G) = \max_{x \in \cC(G)} \nu(x),  \;\;\mbox{and}\;\;\nu_*(G) = \min_{x \in \cC(G)} \nu(x).
\end{equation}

\subsubsection*{Properties of construction sequences}
One may consider three useful kinds of construction sequence \cite{kp}.  Let $G = (V,E)$ be a graph.  Call $x \in \cC(G)$ an {\bf easy} c-sequence if no edge proceeds a vertex.  
There are a total of $p!q!$ easy c-sequences for $G$.   

We say that $x\in \cC(G)$ is a {\bf greedy} c-sequence if it has the property that for all edges $e = uw$ and all vertices $v \notin \{u,w\}$
\begin{equation}
x^{-1}(v) < \max \{x^{-1}(u), x^{-1}(w)\} \;\;\mbox{or}  \;\;x^{-1}(v) > x^{-1}(e).
\end{equation}
A greedy c-sequence is characterized by the condition that once a vertex $w$ has been placed, $w$ must be followed by all the edges available to it before any other vertex can be placed.  That is, if $x_i = w$, then there exists $k \geq 0$ such that for $1 {\leq} j {\leq} k$, $x_{i+j} = u_j w \in E$ with $\{u_{i+1}, \ldots, u_{i+k}\} = \{u \in V: x^{-1}(u) < x^{-1}(w)\}$.  

A c-sequence can be both easy and greedy.  For the star $K_{1,n}$, take all peripheral vertices followed by the hub vertex and then all the edges.

A construction sequence $x \in \cC(G)$ is {\bf nearly connected} if the number of connected components in $G_i$,
$1 \leq i \leq \ell$, never exceeds 2.  Equivalently, a c-sequence (for a connected nontrivial graph) is nearly connected if it contains a unique pair of consecutive vertices. See also \cite{UK}.

A small set of basic lemmas will be used repeatedly in this paper.  Let $\cC^*(G)$ and  $\cC_*(G)$ denote the sets of max and min cost construction sequences for $G$.

A min cost c-sequence must be greedy.
\begin{lemma}
    If $x \in \cC_*(G)$, then $x$ is a greedy sequence.
\end{lemma}

\begin{proof}
If $x \in \cC(G)$ is not greedy, then some vertex $v$ appears between $e=uw$ and the later of its two endpoints, so $\max\{x^{-1}(u), x^{-1}(w)\} < x^{-1}(v) < x^{-1}(e)$. Then $v = x_i$ and $e = x_j$ for $i < j$.  Define $y_k = x_k$ for $k \in [\ell] \setminus \{i,j\}$ and put $y_i = x_j$ and $y_j = x_i$. Then $y \in \cC(G)$ and $\nu(y) < \nu(x)$.
\end{proof}

A max cost c-sequence must be easy.
\begin{lemma}
    If $x \in \cC^*(G)$, then $x$ is an easy sequence.
\end{lemma}

\begin{proof}
If $x$ is not easy, some edge immediately precedes a vertex. Reversing them increases cost. Indeed, if $x_i = e$ and $x_{i+1} = v$, where $e = uw$ and $v \notin \{u,w\}$, then defining $y_i = v$, $y_{i+1} = e$, and $y_j = x_j$ for $j \in [\ell] \setminus \{i, i+1\}$, we have $y \in \cC(G)$ and $\nu(e,y) = 2+\nu(e,x)$; the $y$-cost of edges incident with $v$ will also increase by 2, while other edge costs are unchanged.
\end{proof}

A max cost c-sequence must put vertices in non-increasing order of degree.
\begin{lemma}
If $x \in \cC^*((G)$ and $x|_V \coloneqq \;(v_1, ... v_n)$, then $\deg(v_i) \geq \deg(v_{i+1})$.
\end{lemma}

\begin{proof}
Suppose for some $i < j$ we have $a := \deg(v_i) < \deg(v_j) =: b$.  
Then swapping $v_i$ and $v_j$ increases the construction cost by $(j -i )(b - a) > 0$. 
\end{proof}

The cost of an easy sequence is not affected by permuting the edges.
\begin{lemma}
If $x, y \in \cC(G)$ are easy and $x|_V = y|_V$,
then $\nu(x) = \nu(y)$.
\label{lm:edge-order-invar}
\end{lemma}

\begin{proof}
When we switch two adjacent edges,  costs increase by 2 for the edge being moved further away while the cost of the other decreases by 2. But any edge permutation is a composition of adjacent transpositions.
\end{proof}

The same argument shows that in a greedy sequence, cost is not affected by the order of any set of edges made available by the inclusion of some vertex.

An {\bf isomorphism} of two graphs is a bijection $\varphi$ between vertex sets such that both $\varphi$ and $\varphi^{-1}$ preserve adjacency.  Isomorphic graphs have identical costs.
\begin{lemma}
Let $\varphi: G \to G'$ be an isomorphism.  For $x \in \cC(G)$, we define $y := (y_1, \ldots, y_\ell) = \varphi(x) \in \cC(G')$ by $y_i = \varphi(x_i), 1 \leq i \leq \ell$.  Then $\nu(x) = \nu(y)$.
\end{lemma}
\begin{proof}
Indeed, for every $e \in E_G$, $\nu(e,x) = \nu(\varphi(e),y)$.
\end{proof}
What is the relationship between max and min cost of  graph and subgraph?

\begin{lemma}
If $H \subseteq G$, then $\nu_*(H) \leq \nu_*(G)$ and, if $VH=VG$, $\nu^*(H) \leq \nu^*(G)$.
\label{lm:induced}
\end{lemma}
\begin{proof}
Any construction sequence $x$ for $G$ induces a c-sequence $y := x|_H$ for $H$ and $\nu(y) \leq \nu(x)$.  Let $x_* \in \cC_*(G)$ and $y_* \in \cC_*(H)$. 
Then 
\[\nu_*(H) = \nu(y_*) \leq \nu(x_*|_H) \leq \nu(x_*) = \nu_*(G).\]
If $H$ is a spanning subgraph of $G$, then, for any easy c-sequence $y$ for $H$, there exists a c-sequence $x$ for $G$ with $x|_H = y$ and by Lemma \ref{lm:edge-order-invar} all such $x$ have the same cost.
Let 
$y^* \in \cC^*(H)$ and choose $x \in \cC(G)$ such that $x|_H = y^*$, 
then 
\[\nu^*(G) \geq \nu(x) \geq \nu(x|_H) = \nu(y^*)=\nu^*(H). \]
\end{proof}

\begin{lemma}
Let $G=(V,E)$ be a graph and let $x \in \cC(G)$ be an easy c-sequence for $G$ with $x|_V = x_{(p)} = (v_1, \ldots, v_p)$ and suppose $\deg(v_j) = d_j$.  Then
\[\nu(x) = q(2p + q + 1) - \sum_{j=1}^p j \,d_j \]
\label{lm:cost-of-degseq}
\end{lemma}
\begin{proof}
By Lemma \ref{lm:edge-order-invar}, edge order doesn't affect cost. The vertex $v_j$ at position $j$ in $x$ has distance $(p + \frac{q+1}{2} - j)$ from the average location of an edge.  But $\deg(v_j) = d_j$, so summing the contributions of $v_j$ to the cost of $x$, we get 
\[\nu(x) = \sum_{j=1}^p d_j\Big(p + \frac{q+1}{2} -j\Big) = \Big(p + \frac{q+1}{2}\Big)(2q) - \sum_{j=1}^p j \,d_j\]
\end{proof}

\subsubsection*{Examples}
Let $P_3 = (\{1,2,3\},\{12,23\})$ be the path with 3 vertices.  It is shown in \cite{kp} that $\cC(P_3)=16$.
Then $x = (1, 2, 3, 12, 23) \in \cC(P_3)$ and $\nu(x) = \nu(12,x) + \nu(23,x) = (2 + 3) + (2 + 3) = 10$. In contrast, for $y = (3, 1, 2, 12, 23)$, $\nu(y) = 3 + 6 = 9$ and for $z = (2, 3, 1, 12, 23)$, $\nu(z) = 4 + 7 = 11$.  So $x, y, z$ are easy sequences.  In fact, $z$ is max cost.  
Greedy sequences for $P_3$ include $x' = (1, 2, 12, 3, 23)$ and $y' = (2, 3, 23, 1, 12)$ with costs $\nu(x') = 7$ and $\nu(y') = 8$; in fact, $x'$ is min cost.

If $K_{1, 3} = (\{0,1,2,3\}, \{01, 02, 03\})$, then $x = (0, 1, 2, 3, 01, 02, 03) \in \cC(K_{1,3})$ and $\nu(x) = 7 + 8 + 9 = 24$. Later, we show that this equals $\nu^*(K_{1, 3})$. 
Another construction sequence for $K_{1, 3}$ is (1, 2, 3, 0, 01, 02, 03), which has a cost of 18 but the nearly connected c-sequence $(1, 0, 01, 2, 02, 3, 03)$ has cost 13 which is the minimum construction cost for the star $K_{1, 3}$.  But for $K_{1, 4} = (\{0,1,2,3,4\}, \{01, 02, 03, 04\})$, the c-sequence $(1, 2, 0, 01, 02,3, 03, 4, 04)$ has min cost but is {\it not} nearly connected as it starts with three isolated vertices.

\section{Maximum costs}


We begin with the maximum cost to construct a regular graph $G=(V,E)$ of order $p=|V|$ and size $q=|E|$.  

\begin{theorem}
If $G=(V,E)$ is regular and $x \in \cC(G)$ is any easy c-sequence, then $\nu(x) = \nu^*(G) =q(p+q)$.
\end{theorem}

\begin{proof}
By Lemma \ref{lm:cost-of-degseq}, it suffices to observe that $d_j = 2q/p$ for $1 \leq j \leq p$. Hence, $\nu(x) =  q(2p + q + 1) - \frac{2q}{p}\frac{p(p+1)}{2} = q(p+q)$. 
 \end{proof}  
 
As a consequence, for $C_p$ the $p$-cycle, $K_p$ the complete graph, $K_{n,n}$ the complete bipartite graph, and $Q_d$ the $d$-cube, one has
\begin{align}
\nu^*(C_p)\; &= \;2p^2\\
\nu^*(K_p)\; &= \;\frac{p^4 - p^2}{4}\\
\nu^*(K_{n,n})\; &=\; 2n^3 + n^4\\
\nu^*(Q_d)\; &= \;(2d + d^2) 4^{d-1}.
\label{eq:cube}
\end{align}

Now consider irregular graphs, starting with stars and paths.  

\begin{theorem}
The maximum construction cost of a star is given by 

    \begin{equation*}
        \nu^*(K_{1,n})  = \frac{5n^2 + n}{2} \text{ for } n \geq 2.
    \end{equation*}
    \label{th:star-max}
\end{theorem}

\begin{proof}
    Let $0$ denote the hub, write $j$ for any of the endpoint vertices adjacent to $0$ for $j = 1, ..., n$, and let $e_j := 0j$. 
By Lemmas 2, 3, 4, and 5, the c-sequence
    \begin{equation*}
      x = (0, 1, 2, ..., n, e_1, e_2, ..., e_n)
    \end{equation*}
is maximum cost.
    Each of the $n$ edges is $n$ steps away from its nearer endpoint, giving $n^2$. Cost to 0, the further endpoint, is the sum $\sum_{j=1}^n (n+j) = \frac{3n^2 + n}{2}$.  Combining costs gives $\nu^*(K_{1,n}) = \frac{5n^2 + n}{2}$, or use Lemma \ref{lm:cost-of-degseq}.
\end{proof}
\begin{theorem}
    Maximum construction cost of a path with $n$ vertices is given by 

    \begin{equation*}
        \nu^*(P_n) = 2n^2 - 2n  - 1 \text{ for } n \geq 2.
    \end{equation*}
\end{theorem}
\begin{proof}
Let $P_n = (V,E)$ be the path with $V=[n]$ and with $n{-}1$ edges joining consecutive integers.
We suppose that $1$ and $n$ are the two endpoints of the path.
If $x \in \cC(P_n)$ is of maximum cost, then $x = (v_1, \ldots, v_{n-1}, v_n, e_1, \ldots, e_{n-1})$, $\{v_{n-1}, v_n\} = \{1, n\}$, and $\nu(x) = \nu^*(P_n)$. 
By Lemma 4, edge order is irrelevant and, as in the proof of Theorem 1, vertex order of $\{v_1, \dots, v_{n-2}\}$
is also arbitrary.  Let
$e_1 = [2,3], e_2 = [3,4], \ldots, e_{n-3} = [n-2,n-1]$, $e_{n-2} = [1,2]$, $e_{n-1} = [n-1,n]$.  Now $\nu(e_j,x) = 2n-1$ for $1 \leq j \leq n-3$ and $\nu(e_{n-2},x) + \nu(e_{n-1},x) = 5n-4$, so we have 
$\nu(x) = \nu^*(P_n) = (n-3)(2n-1) + 5n-4 = 2n^2-2n-1$.
\end{proof}
Thus, maximum costs for the two extremal trees grow quadratically with the number of vertices and their ratio approaches the limit $5/4$.  The following theorem says that max construction costs for {\it all} trees grow at quadratic rates.


\begin{theorem}
If $T$ is a tree with $n$ vertices, then $\nu^*(P_n) \leq \nu^*(T) \leq \nu^*(K_{1,n-1})$.
\end{theorem}
\begin{proof}
By Lemma \ref{lm:cost-of-degseq}, the {\it least} max cost for a tree with $n$ vertices occurs when $\sum_{j=1}^n j d_j$ is {\it largest} corresponding to degree sequence $(2, \ldots, 2, 1, 1)$ for paths and similarly the largest max cost occurs for the degree sequence $(n, 1, \ldots, 1)$ for stars.
\end{proof}


\medskip

The disjoint union of regular graphs of the same degree is again regular, so $\nu^*(G_1 \sqcup G_2) = (q_1 + q_2)(p_1 + p_2 + q_1 + q_2) = \nu(G_1) + \nu(G_2) + q_1(p_2 + q_2) + q_2(p_1 + q_1)$.
The next theorem describes the case for regular graphs of distinct degrees.
\begin{theorem}
Let $G_1$ be $r_1$-regular and let $G_2$ be $r_2$-regular with  $r_1 > r_2$. Then

    \begin{equation*}
        \nu^*(G_1 \sqcup G_2) = \nu^*(G_1) + \nu^*(G_2) + 2q_1(p_2+q_2),
    \end{equation*}
where $p_i {:=} |V_i|{=} |V(G_i)|,\;\; q_i{:=}|E_i|{=}|E(G_i)|$ for $i=1,2$, and $\sqcup$ is disjoint union.    
\end{theorem}

\begin{proof}
By Lemmas 2, 3 and 4, if $r_1 >  r_2$, then there exists $x \in \cC^*(G_1 \sqcup G_2)$ that is a 4-fold concatenation, $x = \rho_1 \,\odot\,\rho_2 \,\odot\, \sigma_2 \,\odot \, \sigma_1$ for permutations $\rho_1$ and $\rho_2$ of $V_1$, $V_2$, and permutations $\sigma_1$, $\sigma_2$  of $E_1$ and $E_2$. Hence, for all $e \in E_2$, $\nu(e,x) = \nu(e,y)$, where $y = \rho_2 \,\odot\,\sigma_2$.  However, for each edge $e'$ in $\sigma_1$, we have
\[\nu(e', x) = \nu(e', \rho_1 \odot \sigma_1) + 2(p_2+q_2).\]
Thus, $\nu^*(G_1 \sqcup G_2) = \nu^*(G_1) + \nu^*(G_2) + 2q_1(p_2+q_2)$. 
\end{proof}


Define for $n\geq 1$ (cf. \cite{sie}) the {\bf symmetric double star} $D_n$ as the graph formed by the disjoint union of two copies of $K_{1,n}$ and an edge joining their hubs. \\

\begin{figure}[ht!]
\centering
\includegraphics[width=50mm,scale=0.5]{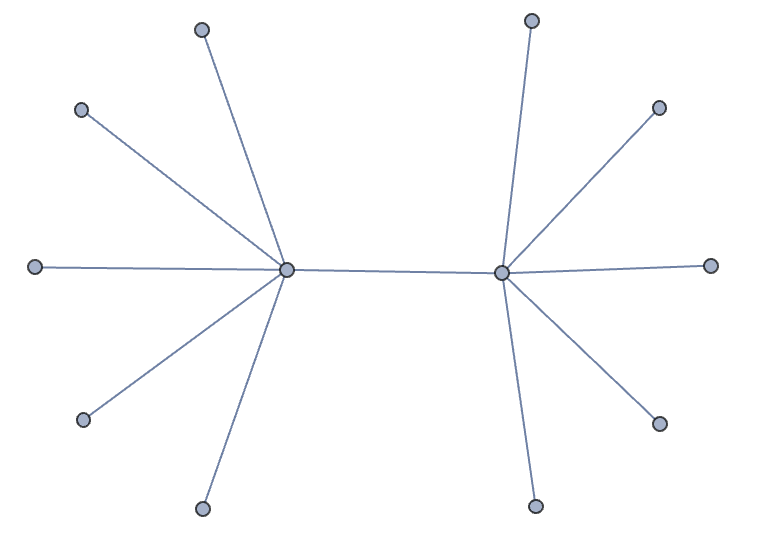}
\caption{The double star $D_5$}
\end{figure}

\begin{theorem}
For $n \geq 1$, $\nu^*(D_n) = 10n^2 + 10n + 3$.
%
\end{theorem}

\begin{proof}
Let $x^* \in \cC(D_n)$ be max cost, $\nu(x) = \nu^*(D_n)$.
Label the two hubs of $D_n$ by $0$ and $0'$, respectively.   Write $j$ for any of the endpoint vertices adjacent to $0$, $j = 1, ..., n$, and put $e_j:= 0j$. Write $j'$ for any of the endpoint vertices adjacent to $0'$ for $j = 1, ..., n$, and let $e_j' := 0'j'$. 
By Lemmas 2, 3, and 4, the two hubs must go first in the maximum construction cost sequence and, by Lemma 5,
we can assume without loss of generality that 
   \begin{equation*}
       x^* =  (0, 0', 1, 2, ..., n, 1', 2', ..., n', e_1, e_2, ..., e_n, e_1', e_2', ... e_n', 00')
    \end{equation*}

    The cost of $e_1$ through $e_n$ is given by the sum of the distance to the closer vertex of distance $2n$ and the distance to the further vertex of distance $2n+1$ plus its index. The sum of cost from $e_1$ through $e_n$ equals $\sum_{j=1}^n 2n+2n+1 + j = n(4n+1) + n(n+1)/2 = 4n^2 + n + (n^2+n)/2$. The cost of $e_1'$ through $e_n'$ is given by the sum of the distance to the closer vertex of distance $2n$ and the distance to the further vertex of distance $3n$ plus its index. The sum of the cost from $e_1'$ through $e_n'$ equals $\sum_{j=1}^n 2n + 3n + j = n(5n) + n(n+1)/2 = 5n^2 + (n^2+n)/2$. Finally, the cost of the edge $e_{0,0'}$ equals $4n+1 + 4n+2$. Therefore, the maximum cost is given by $4n^2 + n + (n^2+n)/2 + 5n^2 + (n^2 + n)/2 + 4n+1 + 4n+2 = 10n^2 + 10n + 3$.
\end{proof}

Let $W_n:=C_n * K_1$ be the $n{+}1$-vertex {\bf wheel graph}, where $*$ means join.

\begin{theorem}
    The maximum construction cost for $W_n$ is given by 

    \begin{equation*}
        \nu^*(W_n)  = \frac{13n^2 + n}{2} \text{ for } n \geq 4.
    \end{equation*}
\end{theorem}

\begin{proof}
Let $y \in \cC^*(W_n)$.  The wheel is a star with peripheral vertices put into a cycle with $n$ new edges.  Using $0$ for the hub, $1, \dots, n$ for the peripheral vertices, and edges $e_j := 0j$, $1 \leq j \leq n$ for the star, and $f_j := [j,j{+}1]$, $1 \leq j \leq n-1$, $f_n := [1, n]$ for the new edges, we can rearrange $y$ to another max cost c-sequence for the wheel, 
$y' = (0, 1, \ldots, n, e_1, \ldots, e_n, f_1, \ldots, f_n)$.
Thus, $\nu(y') = \nu^*(K_{1,n}) + (n-1)(4n-1) + 5n-1 = \frac{5n^2+n}{2} + 4n^2$, as required.
\end{proof}

\section{Minimum cost c-sequences}

Determining {\it minimum} cost of a construction sequence for $G$ seems to be more difficult than the max cost problem. We can solve it for complete graphs.

\begin{theorem}
For $n \geq 1$, $\nu_*(K_n) = \frac{(n-1)n(n+1)(n+4)}{12}$.
\end{theorem}
\begin{proof}
A min cost sequence must be greedy by Lemma 1 and by Lemma 5, any greedy sequence is minimum cost.  Hence, one can choose $x \in \cC_*(K_n)$ of the form
\[x = (1,2,\overline{12},3,\overline{13},\overline{23},4,\overline{14},\overline{24},\overline{34},5,\overline{15} \dots)\]
where vertices are denoted $1, 2, \ldots, n$ and edge $ij$ is written $\overline{ij}$ for clarity.  

It will suffice to show that for each integer $k$, $2 \leq k \leq n$, one has the following sum of edge costs
\begin{equation}
\sum_{j=1}^{k-1} \nu(\overline{jk},x) = k(k-1)(2k+5)/6
\label{eq:cost-of-star}
\end{equation}
and we denote this cubic polynomial by $X_k$.
Summing $X_k$ for $2 \leq k \leq n$, e.g., using the 
computational tool
Wolfram Alpha \cite{wolfram}, one gets $\nu_*(K_n)$ as asserted.

Now we prove (\ref{eq:cost-of-star}). First note that for $j \in [n]$, one has $x_i = j$ exactly when $i=1+{j \choose 2}$. Further, for $2 \leq k \leq n$, the $k{-}1$ edges $\overline{jk}$, $1 \leq j \leq k{-}1$, immediately follow $k$ so if $x_i = j$, then $x_{i+j} = \overline{jk}$. Thus, for $1 \leq j \leq k{-}1$, the cost of edge $\overline{jk}$ is the sum of its delays from $j$ and $k$; that is, 
\[
\nu(\overline{jk},x) = \Big({k \choose 2} + 1 + j - {j \choose 2} - 1 \Big) + \Big(j\Big),
\]
where the last summand $j$ is the distance from $k$ to $\overline{jk}$ for $j < k$ in $x$.   Hence, 
\[\sum_{j=1}^{k-1} \nu(\overline{jk},x) = \frac{k^3 - k^2}{2} + k^2 - k - {k+1 \choose 3} = \frac{(k^2 - k)(3k+6-k - 1)}{6}.\]
Therefore, (\ref{eq:cost-of-star}) is verified and the theorem follows. 
\end{proof}
An interesting consequence is that the ratio of max to min cost for the complete graph is asymptotically equal to 3. For paths, however, the ratio of max to min cost grows linearly with the number of vertices as we now show.
\begin{theorem}
For $n \geq 3$, $\nu_*(P_n) \leq 4n - 5$.
\label{th:max-path}
\end{theorem}

\begin{proof}
Stanley \cite{stan-book} defines a {\bf fence} as the poset determined by the vertices and edges of $P_n$ under inclusion of a vertex in the edge of which it is an endpoint.
Let $j$, $1 \leq j \leq n$ denote the vertices and $\overline{j} := [j, j{+}1]$ , $1 \leq j \leq n-1$ the edges of $P_n$.
In the natural layout of the Hasse diagram of a fence poset, the vertices are listed from left to right in the natural (path) order along level one at positions $1, 3, \ldots, 2n+1$ while  the edges $\overline{1}, \ldots, \overline{n{-}1}$ are listed left to right order at positions $2, \ldots, n{-}1$ along level two with all inclusion arrows going upward.  This gives a fence-like visual pattern.

To linearize the poset with as small as possible a total delay, move each edge one step forward in the left-to-right ordering of the elements of the path.  Instead of the sequence $1, \overline{1}, 2, \overline{2}, \ldots, n{-}1,\overline{n{-}1}, n $, we have $1, 2, \overline{1}, 3, \overline{2}, \ldots, n, \overline{n{-}1}$.  The first edge costs 3 but subsequent edges cost 4. \end{proof}

If the min cost sequence $x$ must be nearly connected, then it must move left to right (or right to left) as above and equality is achieved.  However, for stars, which we examine next, min cost sequences need not be nearly connected.
But this does not seem to be possible for a path, where we claim that a minimum cost c-sequence $x$ must be nearly connected.  

In a minimum c-sequence $x$ for the path $P_n$, let $x_1 = v$.  Either $x_2 = w$ and $x_3 = vw \in E$ or else $x_2 = u$ and $vu \notin E$.  In the second alternative, whenever we place $w$, it will be followed by the edge $vw$.  The cost involved in the additional delay between $v$ and $vw$, when $w$ does not immediately follow $v$ in $x$, is larger than the optimal saving in delayed edges, which is at most 1.  Indeed, an edge has cost 3 if and only if it appears as an isolated component of the evolving path. But the standard order costs 4 per edge. 
If isolated edges occur before $w$, the additional cost for $vw$ is still higher. Hence, $x$ is nearly connected.

In the construction process, it is possible to initially keep cost-per-edge at 3 by taking $\lfloor n/3 \rfloor$ edges such that no two of these edges have adjacent endpoints, but such c-sequences have greater total cost.



\begin{theorem}
    For $n \geq 3$, the minimum construction cost of a cycle is given by 

    \begin{equation*}
        \nu_*(C_n) = 6n - 4.
    \end{equation*}
\end{theorem}

\begin{proof}
We use the notation of the proof of Theorem \ref{th:max-path}.  Let $x$ be a min cost sequence for $C_n$ and let $y = x_{(2n{-}1)}$ be obtained from $x$ by deleting the last edge $e=uw$.  The two vertices $u$ and $w$ must be of degree 1 in a path subgraph $P_n$ of $C_n$ and $y = x|_{P_n}$.  But only the first and last vertex in $y$ are of degree 1 in $P_n$, so the cost of $x_\ell$ is $2n+1$.  
By the proof of Lemma \ref{lm:induced}, \[4n - 5 = \nu_*(P_n) \leq \nu(x|_{P_n}) = \nu(x) - \nu(x_\ell, x)\] 
so $\nu_*(C_n) \geq 4n-5 + 2n+1\geq 6n - 4$;
the ordering for the path gives equality.
\end{proof}

We now consider the star $K_{1,n} = \Big(\{0\}\cup[n], \{\overline{j}: j \in [n]\} \Big)$, where $\overline{j} = 0j$.  By symmetry of the graph, the order of all vertices in $[n]$ is arbitrary, as is the order of the available edges.  The distinguishing feature in a minimum c-sequence for the star is the number $b$ of vertices that come before $0$.  Everything else will be determined.  No edge can be placed before $0$.  If $0$ is placed first, it is followed by another vertex $j$ and then, by the greediness of $x$, by $\overline{j}$ with cost 3.  But the same is achieved by $(j, 0 , \overline{j})$ at equal cost and all later edges $\overline{i}$ are cheaper.  The same argument shows that the initial block of all vertices must have $0$ last.  

Suppose $b=1$. This extreme value requires $x$ to be nearly connected and corresponds to the ``straightforward'' c-sequence
\[x' = (1,0,\overline{1},2,\overline{2}, \ldots, n, \overline{n})\]
which has a cost of 
\[\nu(x') = 1 + \sum_{j=1}^{n} 2j = 1+n(n+1) = 1+ n^2 + n.\]


In contrast, suppose $n = 2r$ and $b=r$ and let 
\[x = (1,2, \ldots, r,0,\overline{1},\overline{2}, \dots,  \overline{r},r{+}1, \overline{r{+}1}, \ldots, 2r, \overline{2r}).\]
The cost of the first $r$ edges is 
\[ \sum_{j=1}^r (r+1+j) =3 {r+1 \choose 2}\] and the cost of the second $r$ edges is $\sum_{j=1}^r (r+1+2j) =4{r+1 \choose 2}$.  Hence, $\nu(x') \sim (7/2)(n^2/4) = \frac{7}{8}n^2$.  This can be still further improved.

\begin{theorem}
Let $n \geq 2$. Then $\nu_*(K_{1,n}) = \frac{5n^2 + 11n}{6}$ and $b = \lceil n/3 \rceil$.
\end{theorem}
\begin{proof}
Suppose $x'' \in \cC(K_{1,n})$ be any greedy c-sequence and let $b$ be the unknown number of vertices preceding the hub $0$.  Edges come in two blocks: those immediately following $0$ and then the remaining edges which each appear between two vertices, except for the last edge.  The cost for the first block is $b(b+2) + b(b-1)/2 = (3/2)(b^2+b)$. The cost for the second block of edges is $(n-b)(b+3) + (n-b)(n-b-1) = (n-b)(n+2) = n^2 + 2n - b(n+2)$. Now replace $b$ by an indeterminate $u$ to give the cost as a polynomial function of $u$, 
\[\nu(x'') = (1/2)\Big(3u^2 + 3u - 2u(n+2) + 2n^2 + 4n\Big),\]
with $n$ as a parameter. So $\nu(x'') = (1/2)(3u^2 - u(1+2n) +2n^2+4n)$.  Setting the $u$-derivative of this expression equal to $0$ (ignoring the factor of $1/2$), we get
\begin{equation}
6u = 1+2n
\end{equation}
at the minimum, giving $u = \lceil n/3 \rceil$.

For $n = 3j$ with first $j$ vertices selected, followed by the hub, then the rest greedily, one gets $\frac{5n^2 + 11n}{6}$
as the min cost.
\end{proof}
At the other extreme value of $b=n$, the cost of a c-sequence is $\sim 3n^2/2$.


We have $\nu_*(K_{1,n}) \sim 5n^2/6$ while from Theorem \ref{th:star-max}, $\nu^*(K_{1,n}) \sim 5n^2/2$ so for stars, the asymptotic ratio of max to min cost is again 3, which is also the asymptotic ratio max to min for complete graphs.

One convenient aspect of min cost is that it is additive over disjoint graphs;  if $x$ is min cost for $G$ and $y$ is min cost for $H$, then $x \odot y$ is min cost for $G \sqcup H$.
\section{Discussion}

We have focused on construction costs for simple undirected graphs, but c-sequences and cost functions can be extended to include multigraphs (with multiple parallel edges joining a given pair of vertices) and digraphs, where edges are assigned a fixed direction. One method for adjusting the construction sequence for directed acylic graphs would be to place a restriction that the vertex order of the endpoints of an arc must follow the direction of the arc.

 Construction costs can serve as an estimate for the complexity of graphs and could be applicable in operations research and optimization.  Minimum cost c-sequences are rather like the ``just-in-time'' notion in logistics \cite{jit}.  
 
 Max cost gives a worst-case bound but might be advantageous if the cost of nodes is cheap now while the cost of links will be cheap later.  In an adversarial situation, it might be good to use a random order to minimize predictability - e.g., a max cost sequence has first nodes with the largest degree so an opponent could immediately attack them.
 
In some experiments, we compared construction cost for Breadth-First Search (BFS) and Depth-First Search (DFS) algorithms. While both algorithms have a run-time of $O(V+E),$ for graphs with the same number of vertices and edges, it may be difficult to immediately compare their runtime. However, the maximum construction cost gives us a method to estimate a comparison between the runtime of two graphs with the same or similar number of vertices and edges. This same idea may be able to be used as a heuristic method across different aspects of operations research and serve as a tool to compare the complexity of projects. 

In this paper, we use a cost based on the integer distances between edges and their endpoints in a construction sequence.  It would be possible to take some power of these distances, which could be different for the nearer and further endpoints.  Such a nonlinear cost function would not allow our convenient lemmas to be used.  One might also have some weight associated with each vertex and edge. Adjusted construction costs could be created to tailor the application of cost to construction sequences in different situations. 

Small graphs, sometimes called ``graphlets'', are appearing in neural networks \cite{liu}, \cite{bio-inf}, \cite{xenos}.
One might think of this as a kind of structured finite elements approach.  Tasks arise in determining which graphs appear and in measuring the difference between two large time-varying graphs, the Worldwide Web at Friday noon EST December 6, 2024 or one week later.  Better would be to identify some sort of macrostate in a large network: Time to sell stocks?

Costs of construction sequences can be used to identify the graph.  For max cost, only the degree sequence matters, so cubic trees of equal order will have equal max cost but taking graph powers extinguishes this ambiguity in most cases.  For example, the two trees $T_1, T_2$ in Harary \cite[p 62, Fig. 6.3]{harary} have the same max cost but $\nu^*(T_1^2) \neq \nu^*(T_2^2)$, where $G^2$ is the supergraph obtained by giving $G$ new edges that join 
vertex-pairs at distance 2
\cite[p 14]{harary}.
Minimum cost is more sensitive.  Already $\nu_*(T_1) \neq \nu_*(T_2)$ even though $T_1$ and $T_2$ have the same degree sequence.  

Lemma \ref{lm:cost-of-degseq} gives a simple formula for max cost so combining it with graph powers might provide better discrimination for a given computational cost than finding minimum cost.


\section*{Appendix}

We evaluate maximum cost to build four additional graph families. Arguments are similar to our earlier results.
The {\bf suspension} $\Sigma G$ of a graph $G$ is the join 
$\Sigma G := G * \overline{K_2}$ of $G$ with a set of two isolated vertices \cite{opx};
e.g., $\Sigma K_n = K_{n+2}$.

\begin{figure}[ht!]
\centering
\includegraphics[width=50mm,scale=0.5]{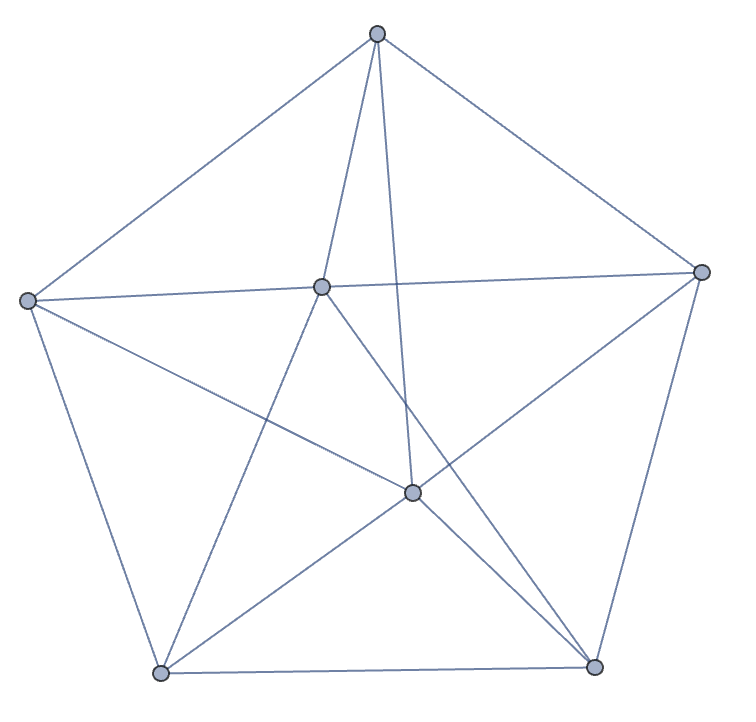}
\caption{The suspension $\Sigma C_5$ of a 5-cycle}
\end{figure}

\begin{theorem}
For $n \geq 4$, $\nu^*(\Sigma C_{n}) = 13n^2 + 2n$.
\end{theorem}
The {\bf 2-wheel with an axle} $W_{n, n}$ is the graph consisting of the union of two copies of the (ordinary) wheel $W_{n}$ together with an edge joining their hubs. \\
\begin{figure}[ht!]
\centering
\includegraphics[width=100mm,scale=0.5]{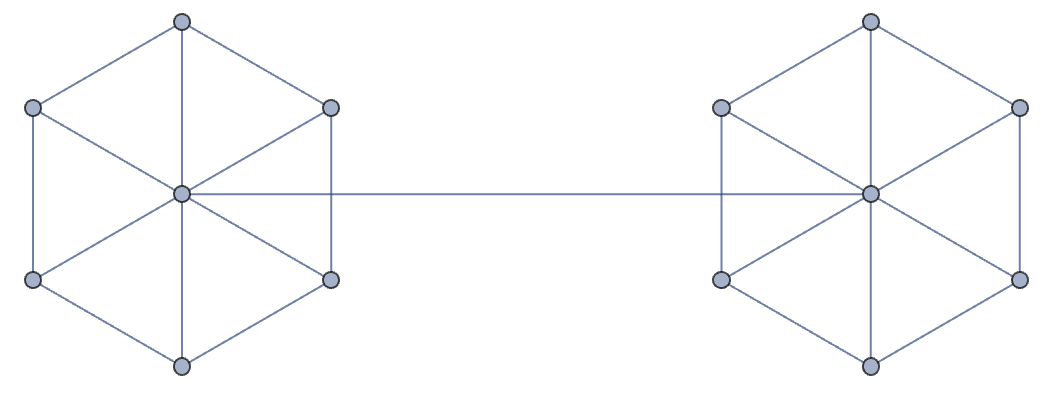}
\caption{$W_{6, 6}$}
\end{figure}

\begin{theorem}
For $n \geq 4$, $\nu^*(W_{n,n})  = 26(n-1)^2 + 14(n-1) + 3$.
%
    \end{theorem}


The {\bf gear} graph $G_n$ is $W_n = K_1 * C_n$ with each edge of $C_n$ subdivided \cite{gj}.\\



\begin{theorem}
For $n \geq 4$, $v^*(G_{n})  = 16n^2 + \frac{n^2+n}{2}$.
%
    \end{theorem}




For $n \geq 1$, the wedge (or ``bouquet'') of triangles $F_n$ is created by taking $n$ copies of $C_3$ with a single common point \cite{gj}. This graph is also known as the $n$-th {\it friendship} graph as it solves the Friendship Problem: characterize the graphs in which each distinct pair of vertices has a unique common neighbor \cite{ers}.\\

\begin{theorem}
For $n \geq 1$, $\nu^*(F_{n}) = 17n^2 + n \text{ for } n \geq 1$.
%
%
\end{theorem}

\begin{proof} A max cost construction sequence for $F_n$ follows that of the star $K_{1, 2n}$, with $n$ additional edges with total cost $\sum_{i=1}^n (8j + 1) - 2i = 7 n^2$ and adding this to $\nu^*(K_{1,2n}) = 10n^2 + n$ gives the result.
 \end{proof}

\end{document}